\documentclass[12pt]{article}
\usepackage{latexsym}
\usepackage{amsfonts}
\makeatletter

\def\@footnotetext#1{\insert\footins{%
\footnotesize
    \interlinepenalty\interfootnotelinepenalty
    \splittopskip\footnotesep
    \splitmaxdepth \dp\strutbox \floatingpenalty \@MM
    \hsize\columnwidth \@parboxrestore
   \edef\@currentlabel{\csname p@footnote\endcsname\@thefnmark}\@makefntext
    {\rule{\z@}{\footnotesep}\ignorespaces
      #1\strut}}}
\def\abstract{\small\quotation{\hskip-\parindent\sc Abstract.}}

\def\classification{\@ifnextchar [{\@xfootnotenext}%
   {\begingroup\let\protect\noexpand
      \xdef\@thefnmark{}\endgroup
    \@footnotetext}}

\title {}

\begin{document}
\classification {{\it 1991 Mathematics
Subject Classification:} Primary 20F28,  secondary 20E36, 57M05.}

\begin{center}
{\bf \large AUTOMORPHISMS  OF  ONE-RELATOR  GROUPS} 
\bigskip

{\bf  
  Vladimir Shpilrain}

\end{center} 
\medskip

\begin{abstract}
\noindent It is a well-known fact that every group  $G$ has a 
presentation of the form $G = F/R$, where $F$ is a free group 
and $R$ ~the kernel of the natural epimorphism from $F$ onto $G$. 
 Driven by the desire to obtain a similar presentation of the 
group of automorphisms $Aut(G)$, we can consider the subgroup 
 $Stab(R) \subseteq Aut(F)$ of those automorphisms of $F$ 
that stabilize $R$, and try to figure out if the natural
homomorphism $~Stab(R) \to Aut(G)$ is {\it onto}, and if it is, 
to determine its kernel.

 Both parts of this task are usually quite hard. The former 
part received considerable attention in the past, whereas the
latter, more difficult, part (determining the kernel) seemed 
unapproachable. Here we approach 
 this  problem for a class of one-relator groups with a special kind of 
 small cancellation condition. 
 Then, we address a  somewhat 
easier case of 2-generator (not necessarily one-relator)  
groups, and determine the kernel 
of the above mentioned   homomorphism  for a rather general class of those 
groups.

\end{abstract}

\date{}

\bigskip

\noindent {\bf 1. Introduction }
\bigskip

 \indent     Let  $F = F_n$   be the free group of a finite rank  
$n \ge 2$  with a  set  
$X = {\{}x_i {\}}, 1 \le i \le n$,  of free generators. Let $R$ be 
a normal subgroup of  $F$, and $Stab(R) \subseteq Aut(F)$ the 
group of those automorphisms $\varphi$  of $F$ that stabilize $R$, 
i.e., $\varphi(R) = R$ (this does not necessarily mean that $\varphi$ 
fixes every element of $R$).  

 Then there is the natural
homomorphism $\rho: Stab(R) \to Aut(G)$, where $G = F/R$ ~(when we
say  {\it the} natural homomorphism, it means every element 
$s \in F$ maps onto the coset $sR$, and this extends to the 
mapping between groups of automorphisms). 
In some important cases, this homomorphism $\rho$ is known to be 
{\it onto}. This was established by Nielsen 
for   surface groups; by Zieschang for somewhat more general
one-relator and Fuchsian 
groups (see \cite{Z} and \cite{Zbook}), and, more
recently, a number of results in this direction have been
published; we only mention a couple of them 
which are in line with the subject
of the present paper. Lustig, Moriah and Rosenberger \cite{LMR}  
completely characterized Fuchsian groups with the  property above.  
Bachmuth, Formanek and Mochizuki \cite{BFM} established this 
property for a  class of 2-generator groups. On the
other hand, it is known by now  that for many groups $G$, the 
homomorphism $\rho$ ~is {\it not} onto -- see \cite{GupSh} or
\cite{Ros2} for a survey. 

 When   $\rho: Stab(R) \to Aut(G)$ happens to be 
onto, it is natural to try to determine its kernel, since that 
would give a presentation of the group $Aut(G)$ as
$Stab(R)/Ker(\rho)$, ~and this, in turn, 
 could give a presentation of $Aut(G)$ 
by generators and defining relations (at least in the case where 
 $G$ is finitely presented) if one uses results of McCool
\cite{McCool} on presentation of  $Stab(R)$. 
\smallskip

 It is clear that the group $Inn_{_{\footnotesize R}}$ of inner automorphisms 
of $F$   induced by elements of  $R$, is always contained in $Ker(\rho)$. 
In many cases, this is the whole $Ker(\rho)$ (see Theorem 1.1 below). 
However, in several important situations  $Ker(\rho)$ happens to be bigger 
than $Inn_{_{\footnotesize R}}$. This is the case, in particular, 
with surface groups: 
\medskip

\noindent {\bf Example 1.(a)} Let $~r = x_1^2 x_2^2 ... x_n^2$, 
~and $~\varphi : ~x_1 \to x_1 \cdot r^{-x_1^{-2}} = 
x_1^{-1}  \cdot x_n^{-2} \cdot ... \cdot x_2^{-2}; ~x_i \to x_i, 
~i \ne 1$ (our notation is $~x^y = yxy^{-1}$). 
 Then $\varphi$ is obviously a non-inner automorphism of the 
free group $F_n$, but it induces the identical  automorphism of the 
group $F_n/R$. 
\smallskip

\noindent {\bf (b)} A similar automorphism works in the orientable
case,   where $~r = [x_1, x_2] \cdot [x_3, x_4] \cdot ... \cdot
[x_{2m-1}, x_{2m}]$ ~(we assume $n = 2m \ge 4$, and our commutator 
notation is  $[x,y] = x^{-1}y^{-1} x y$):  

\noindent $\varphi_o : ~x_1 \to x_1 \cdot r = x_2^{-1}x_1x_2 \cdot
[x_3, x_4] \cdot ...  \cdot [x_{2m-1}, x_{2m}]; ~x_i \to x_i, 
~i \ne 1$. 

\medskip

\noindent {\bf Example 2.(a)} A little more sophisticated 
example is provided (in the non-orientable case) 
by the following automorphism (here $~r = x_1^2
x_2^2 ... x_n^2$): 

\noindent $~\psi : ~x_1 \to r^{-x_1^{-2}} \cdot x_1 \cdot 
r^{-x_n^{-2}} \cdot r^{x_1^{-2}}; ~x_i \to x_i, ~i \ne 1,n; 
~x_n \to r^{-x_1^{-2}} \cdot x_n \cdot r^{x_1^{-2}}$. 
\smallskip

\noindent {\bf (b)} Based on a similar idea, we get 
the following automorphism in the orientable case: 

\noindent $~\psi_o : ~x_1 \to r^{[x_2,x_1]} \cdot x_1; 
~x_2 \to x_2^{r^{[x_2,x_1]}}; ~x_i \to x_i, ~i \ne 1,2$. 

\medskip 

 It is not quite obvious that the 
automorphisms from  Examples 1 and 2 
belong to $Stab(R)$; it is however obvious that 
we have $\varphi(R) \subseteq R$ for any of them; now $\varphi(R)= R$
 ~follows from the fact that surface groups are hopfian (this was
established by Hopf himself).  
 \smallskip 

 These are not the only examples of non-inner automorphisms of  a 
free group $F$ that induce the identical  automorphism of the 
corresponding surface group $F/R$. In fact, the group 
$Ker(\rho)/Inn_R$ ~appears to be rather big; at least, it is
non-periodic and non-abelian (the automorphisms $~\varphi$ 
and $~\psi$ from  Examples 1(a) and 2(a) do not commute modulo 
$Inn_{_{\footnotesize R}}$). The problem of an actual description
of the group  $Ker(\rho)$ for any of the surface groups (of rank $ >
2$)  remains  open (to the best of my knowledge).  We note that 
 Dicks and  Formanek \cite{WF} recently described (in particular) 
the kernel of a 
 natural homomorphism (``collapse") from the automorphism 
group of the fundamental group of an orientable surface with 
one puncture onto the corresponding (standard) surface group. 
 \smallskip 

 Here we suggest a possible approach to the problem of describing 
$~Ker(\rho)$,  which works 
for many one-relator groups that satisfy a strong type of 
small cancellation condition, in particular for groups that 
 are in some sense close to  surface groups. 
But, as it happens also in some other 
instances, the situation with 
surface groups appears to be on the border between 
difficult and impossible, so it is not clear at the moment 
if this approach can be pushed through  for surface groups as well. 
 \smallskip 

 Before we give the statement of our first result, we need to say 
a couple of words about the Whitehead  graph $Wh(u)$ of a free
group  element $~u$. The vertices of this graph correspond to the
elements of the generating set $X$ and their inverses. If the 
word $~u$ has a subword  $x_i x_j$, then there is an edge in 
$Wh(u)$ that connects the vertex $x_i$ to the vertex $x_j^{-1}$;
~if $~u$ has a subword  $x_i x_j^{-1}$, then there is an edge 
that connects $x_i$ to $x_j$, etc. ~We note that usually, there 
is one more edge (the external edge) included in the definition 
of the Whitehead  graph: this is the edge that connects the vertex 
corresponding to the last letter of $~u$, to the vertex 
corresponding to the inverse of the first letter. We shall not 
include  the external edge in $Wh(u)$; instead, we shall consider 
the Whitehead  graph $Wh(\overline{u})$ of a {\it cyclic} word 
$~\overline{u}$ when  necessary, in which case, of course, the
external edge is 
 included. 

\medskip 

\noindent {\bf Theorem 1.1.} Let  $G=F_n/R, ~n \ge 3,$ be a
one-relator group  with  a  relator $~r~$ and with the following
property: $G$ satisfies  a small cancellation condition
$C'(\lambda), ~\lambda \le 1/6$, and the Whitehead  graph of any
subword  of length $ \ge (1 - 3 \lambda)|r|$ ~of the word 
 $~r~$ or any of its cyclic permutations, is 2-connected (i.e., is 
connected and  does not have a cut vertex). Let  $\rho: Stab(R) \to
Aut(G)$ be the natural homomorphism.  Then $Ker(\rho) = 
Inn_{_{\footnotesize R}}$.  

\medskip 

 If  $n = 2$,  there might be  some additional automorphisms in 
$Ker(\rho)$, namely, inner automorphisms induced by elements $u \in
F$   such 
  that   $~uR~$ is central in $G=F/R$. We treat 2-generator groups 
separately, in Theorem 1.3. 
\medskip

 We note that it is easy and straightforward to check if a given 
one-relator group satisfies the conditions of Theorem 1.1. 

 Although surface groups do not satisfy all those  conditions, there are ``similar" groups that do; those are, for
example, one-relator groups with the relator of the form 
$~(x_1^2 x_2^2 ... x_n^2)^p, ~p \ge 2$, ~or $~([x_1,
x_2] \cdot  [x_3, x_4] \cdot ... \cdot [x_{2m-1}, x_{2m}])^p, ~p \ge
2$. (In general, one-relator groups with a ``very long" relator
tend to satisfy those conditions). 

 Note also that the stabilizer of a (cyclic) word 
$~(x_1^2 x_2^2 ... x_n^2)^p,$ ~or $~([x_1,
x_2] \cdot  [x_3, x_4] \cdot ... \cdot [x_{2m-1}, x_{2m}])^p$
 ~is the same as that  of $~x_1^2 x_2^2 ... x_n^2$, or
$~[x_1, x_2] \cdot  [x_3, x_4] \cdot ... \cdot [x_{2m-1},
x_{2m}]$,  ~respectively, and  that the homomorphism 
$\rho: Stab(R) \to Aut(G)$ is {\it onto} for any of the
corresponding  one-relator groups $G = F/R$ ~(see \cite{Ros1} or 
\cite{Ros2}). Therefore, we have:  
\medskip

\noindent {\bf Corollary 1.2.} Let  $G=F_n/R, ~n \ge 3,$  be a
one-relator group  with the relator of the form $~(x_1^2 x_2^2 ...
x_n^2)^p, ~p \ge 2$, ~or $~([x_1, x_2] \cdot  [x_3, x_4] \cdot ...
\cdot [x_{2m-1}, x_{2m}])^p, ~p \ge 2$. Then $Aut(G) =
Stab(R)/Inn_{_{\footnotesize R}}$. 

\medskip

  Theorem 1.1 is proved  the following  way: first we apply small
cancellation theory to make sure that we have a sufficiently large 
fragment of $r^{\pm 1}$ (or some of its conjugates) in every element
of $R$, and  then use this large fragment to show that 
the Whitehead graph of a cyclically reduced 
 element of the form $~x_i \cdot s, ~s \in R$, cannot be the
Whitehead graph of a primitive element  of a free group,
since those large fragments appear to be ``primitivity-blocking" 
because the Whitehead graph of such a fragment does not have a cut 
vertex.  This latter  observation  is essentially 
 due to E.Turner (informal communication). 
\smallskip 

  One more remark about automorphisms of surface groups is 
in order. The group of {\it outer} automorphisms of a surface 
group is known to be isomorphic to the {\it mapping class group } 
of the corresponding surface, and this latter group was studied 
extensively by a number of people. In
particular, Lickorish, Birman and others obtained several 
 different presentations of mapping class groups based on 
various geometric ideas (see \cite{Birman}). Recently, a very 
``economical" presentation was  found by Wajnryb \cite{Wajnryb}.
 These presentations however do not help much in obtaining a 
presentation of the whole group of automorphisms of a surface 
group. 
 \medskip

 Finally, we consider the problem of determining $Ker(\rho)$ for 
2-generator (not necessarily one-relator) groups. These groups are
easier to handle because  of a very convenient criterion of
primitivity for an element of  the free group of rank 2, which can
be found in \cite{CMZ} and  (somewhat disguised) in \cite{OZ}. Based
on this criterion, we prove: 
\medskip 

\noindent {\bf  Theorem 1.3.} Let  $R \subseteq [F_2, F_2]$. If 
$\varphi \in Ker(\rho)$  (where, as before, $\rho: Stab(R) \to  Aut(G)$
is the natural  homomorphism), then  $\varphi$ is an inner 
automorphism induced by an element $s \in F_2$   such 
  that   $~sR~$ is central in $F_2/R$. 

\medskip

 In particular, if a group $G=F_2/R$ ~has trivial centre, then
 $Ker(\rho) = Inn_{_{\footnotesize R}}$. 
\smallskip 
 
 We emphasize once again that 
   when the rank of a group $G $ is bigger than 2, the situation 
becomes much more complicated. In particular, it is not known 
what $Ker(\rho)$  is when $G $ is a free metabelian group of rank 
$> 2$. It is not even known whether or not $Ker(\rho)$  is finitely 
     generated as a normal subgroup of  $Aut(F)=Stab(R)$ in that
case. What is  known is that $\rho$  is {\it onto} for a free
metabelian group  of rank $> 3$ (see \cite{BM}). We draw a special
attention to free  metabelian groups here because the problem of
determining $Ker(\rho)$  for those groups is closely  related to a
notorious problem of  combinatorial group theory  and algebraic
topology -- to the problem of  the Gassner representation of a pure
braid  group being  faithful (see \cite{Birman}, Section 3.3).\\

\noindent {\bf 2.  One-relator groups} 
\bigskip

\noindent {\bf Proof of Theorem 1.1.} We start by recalling a 
well-known property of  the Whitehead graph of a free group element.
If an element  is primitive and cyclically reduced, then its
Whitehead graph has a cut vertex, i.e., a vertex that, having been
removed from the graph  together with all 
incident edges, increases the number of connected components of the 
graph. 

 Therefore, if we want to prove that some element of a free group is 
{\it not} primitive, it is sufficient to show that this element has 
a subword whose Whitehead graph is 2-connected. 
\smallskip 

 Now, by way of contradiction, suppose there is an 
automorphism $\varphi$  of the free group $F_n$ that takes 
$x_i$ to $x_i \cdot s_i, ~1 \le i \le n$, ~where $s_i \in R$, 
and at least one of the $s_i$, say, $s_1$, is non-trivial. 

 We need to have some of the elements $x_i \cdot s_i$ 
cyclically reduced to apply the property of  the Whitehead graph 
discussed above. Suppose, for example, that $u_1 = x_1 \cdot s_1$
is cyclically reduced (this means, in particular, that there is 
no  cancellation between $x_1$  and $s_1$).
Then,  $u_1$ being primitive, the Whitehead 
graph $Wh(u_1)$ must have a cut vertex. However, there is a subword
of  $u_1$ whose Whitehead graph has no cut vertex. Indeed, a 
small cancellation condition $C'(\lambda), ~\lambda \le 1/6$
implies that every element of $R$ has a subword whose length is 
more than $~(1 - 3 \lambda)|r|$  and which is a subword of some 
cyclic permutation of   $~r^{\pm 1}$ ~-- see \cite{LS}, Theorem
V.4.4. Now the conditions of Theorem 1.1 imply that the Whitehead 
graph of such a subword has no cut vertex, and therefore neither 
does $Wh(u_1)$, hence a contradiction. 

 If there is a cancellation between $x_1$  and $s_1$, i.e., if 
$s_1$ is of the form  $x_1^{-1}s'_1$ (but not of the form  
$s'_1x_1^{-1}$), then we might lose one letter in our long 
subword described in the previous paragraph; that is why we 
require the condition on the Whitehead graph to be satisfied by 
subwords of length $ \ge (1 - 3 \lambda)|r|$, not just 
$ > (1 - 3 \lambda)|r|$ as it appears in Theorem V.4.4 
of \cite{LS}. 

\smallskip 

Suppose now that all of $x_i \cdot s_i$ are cyclically reducible; 
that means, in particular, that there might be cancellations between 
$x_i$ and $s_i$.  Then start composing our automorphism $\varphi$ 
with   automorphisms $\psi_{i,k}$, where $k$ runs through integers,
$2 \le i \le n$, 
and $\psi_{i,k}$  takes $x_1$ to $x_1 x_i^k$ and fixes other 
generators. Thus, $\psi_{i,k} \circ \varphi$   takes $x_1$ to 
$v_1 = x_1 \cdot s_1 \cdot (x_i \cdot s_i)^k$. 
If  for some $i, k$ the
element $v_1$ is   cyclically reduced  (of course, it is still
primitive),  then we are done since the previous  argument applies.
 
 If  for  all pairs $(i, k)$ the
element $v_1$ is not  cyclically reduced, this can only mean that 
every  $\varphi(x_i) $ has a form $w_i^g = (x_i \cdot y_i)^g$ for
some $g \in F, ~y_i \in R$.  In that case, we see that the
element $x_i \cdot y_i$ is itself primitive (as a conjugate of a 
primitive element), but we have seen that this is only possible 
if $y_i = 1$, i.e., if our automorphism $\varphi$ was just the 
conjugation by $g$. 

 Thus, we have shown so far that  only inner automorphisms of 
a free group might belong to $Ker(\rho)$. Since every $n$-generator 
one-relator group has trivial centre provided $n \ge 3$ (see 
\cite{LS}, Proposition II.5.22), this implies $Ker(\rho) = 
Inn_{_{\footnotesize R}}$, 
 and this 
completes the proof of Theorem 1.1. \\

\noindent {\bf 3. Two-generator groups} 
\bigskip

\noindent {\bf Proof of Theorem 1.3.} We start by recalling a convenient 
necessary condition of primitivity in $F_2$ (see \cite{CMZ}): 
\smallskip 

\noindent -- if $w$ is
a primitive element of $F_2$, then some conjugate of $w$ can be
written in the form $x_1^{k_1} x_2^{l_1} ... x_1^{k_m} x_2^{l_m}$,
so that some of $x_i$ occurs either solely with exponent 1 or 
solely  with exponent $-$1. 
\smallskip 

 To prove Theorem 1.3, it is sufficient to prove the following 
statement:  if a 
cyclically reduced primitive element of $F_2$ has a form  $x_1 \cdot
c,  ~c \in [F_2, F_2]$, ~then $~c = 1$. Indeed, if we prove it, it
will  follow that the only situation where an automorphism of $F_2$
may  induce the identical automorphism of $G$, is where both
generators  $x_1$ and $x_2$ are taken to their conjugates. In that
case, by a  well-known result of Nielsen (see e.g. \cite{LS}, 
Proposition I.4.5), this automorphism  must be inner. The result
follows.

 Now we prove the statement. If $~c \in [F_2, F_2]$, then both $x_1$ 
and $x_2$ occur in $c$ both with positive and negative exponents. 
Therefore, 
the only way for an element of the form $x_1 \cdot c$ to be primitive 
is to have a cancellation of  $x_1$ with the first letter of $~c$, 
which has to be $x_1^{-1}$. Thus, let $c = x_1^{-1} \cdot c_1$, where 
$~c_1$ has only positive occurences of $x_1$. But this is possible only 
if $c_1 = x_2^k x_1 x_2^{-k}$; ~in that case, we have $x_1 \cdot c = 
x_2^k x_1 x_2^{-k}$, ~a cyclically reducible element, contrary to 
our assumption. This completes the proof.\\

\noindent {\bf Acknowledgement}
\smallskip

\indent I am grateful to J.McCool for  useful discussions.

\medskip
\noindent 
 Department of Mathematics, The City  College  of New York, New York, 
NY 10031 
 
\smallskip

\noindent {\it e-mail address\/}: shpil@groups.sci.ccny.cuny.edu

\end{document}